\documentclass{article}

\usepackage{amsmath}
\usepackage{color}
\usepackage[letterpaper,breaklinks,pdftex,bookmarks,plainpages=false,
   colorlinks]{hyperref}

\title{Hypergeometric tail inequalities:\\ending the insanity%
\footnote{First version posted on the Web, March 2009.  Revised for minor
typographic errors, February 2011.  Revised for posting on arXiv,
November 2013.}}
\author{Matthew Skala\\
   \href{mailto:mskala@ansuz.sooke.bc.ca}{mskala@ansuz.sooke.bc.ca}}

\allowdisplaybreaks

\begin{document}

\maketitle

\section{Introduction}

I recently needed to put a tail inequality on an hypergeometric
distribution. This should be an easy thing to do; but I found the available
online sources to be really frustrating.  Everybody uses different notation,
and most people seem to like giving helpful examples in which the word
``success'' is used to describe failure and vice versa, and the whole thing
is likely to drive the reader nuts.  Here, for my own future reference and
for the benefit of anyone trying to do the same thing, is a summary of what
I was able to glean in what I hope will be clearly understandable terms.

In the years since 2009, when I first posted these notes on my Web site,
they have attracted a fair bit of attention and even some citations in
serious academic publications, not all of which spelled my name correctly. 
Thus it seems appropriate to post the notes on arXiv to make future
citations easier, increase my own visibility in academic search engines, and
so on.

I don't claim there's any original math in these notes; this is just a
summary of well-known results; but it cost me a fair bit of annoyance to get
issues like notation straightened out.  If you use these notes, a citation
to this posting on arXiv would be appreciated.

It is assumed that you know about as much as I did about this stuff before I
did the research: namely, you should know enough to know that applying a
tail inequality to an hypergeometric distribution is what you want to do,
even if you have trouble keeping track of the parameters of the distribution
or knowing exactly which tail inequality you want.  You're also expected to
be mentally flexible enough to translate the balls-and-urn description into
whatever your real application is.  I'll spare you the confusing
burned-out-lightbulbs example.  My own actual application had to do with
counting bits in the bitwise AND and OR of random bit strings with known
numbers of $1$ bits.

The articles by Chv\'atal and Hoeffding may be hard to find online,
especially if you don't have academic library
privileges~\cite{Chvatal:Tail,Hoeffding:Probability}.  Contact me by email
if you need help locating them.

\section{Setup and notation}

You've got an urn with $N$ balls in it.  Some of them, namely $M$ of them,
are white.  The rest, namely $N-M$ of the balls, are black.  You're going to
draw out $n$ balls from the urn.  You are drawing them {\em uniformly},
which means that every time you pull out a ball it is equally likely to be
any of the balls in the urn at that moment.  However, you are drawing them
{\em without replacement}, which means that after you've drawn out a ball of
one colour, you've reduced the number of balls of that colour remaining and
so the next one will be a little more likely to be the other colour.  If
instead you threw each ball back in after drawing it, then every draw would
have the same chances, we'd be dealing with the geometric distribution
instead of the hypegeometric distribution, and the math would be a lot
easier.  But this time you're drawing without replacement.

Now, how many white balls are you going to get among the $n$ you draw? 
Let's call this number $i$; the question is what interesting things we can
say about the distribution of the random variable $i$, which is called an
hypergeometric distribution.

The short answer is that you will get about the same fraction of white balls
in your $n$-ball sample as the fraction of white balls among the $N$ that
the urn contained at the start.  That's the expected value of $i$. 
Moreover, you will nearly always get very close to exactly that fraction. 
The distribution has light tails.  It isn't a normal distribution bell curve
(which is approached by a geometric distribution, which in turn is what
you'd get by sampling with replacement) but it does have the same kind of
faster-than-exponential fall-off that you would get from the normal
distribution.  As a result you can put a limit just a little bit above the
expected value of $i$ and say ``$i$ is nearly always below this limit'' or
put another limit just a little below and say ``$i$ is nearly always above
this limit.'' That is what a tail inequality does.

The usual suspects (MathWorld~\cite{MathWorld} and
Wikipedia~\cite{Wikipedia}) and their sources use many different notations. 
I am following the notation for variables
used by Chv\'atal~\cite{Chvatal:Tail}, because
his paper seemed easiest to understand.  If you try to read the encyclopedia
entries, you can try to translate using this table:

\noindent
\begin{tabular}{lccc}
   & Chv\'atal~\cite{Chvatal:Tail} & MathWorld~\cite{MathWorld} &
     Wikipedia~\cite{Wikipedia} \\
   & and these notes & & \\
   balls in urn           & $N$   & $n+m$ & $N$ \\
   balls that count       & $M$   & $n$   & $m$ \\
   balls that don't count & $N-M$ & $m$   & $N-m$ \\
   balls you draw         & $n$   & $N$   & $n$ \\
   drawn balls that count & $i$   & $i$   & $k$
\end{tabular}

\section{The distribution}

What's the chance of getting exactly $i$ white balls?  For that we want the
probability distribution function; Chv\'atal doesn't give a notation for it
but I am using one based on his notation for the cumulative distribution
function:
\begin{equation}
   h(M,N,n,i) = \binom{M}{i}\binom{N-M}{n-i}\left/\binom{N}{n}\right.
\end{equation}

That follows from simple counting: how many ways can we draw out $n$ balls
including exactly $i$ of the $M$ white balls, compared to the number of ways
we can draw out $n$ balls without caring about how many of them are white? 
The answer is that we must choose $i$ of the $M$ white balls to draw, hence
the factor of $\binom{M}{i}$, and $n-i$ of the $N-M$ black balls, hence
$\binom{N-M}{n-i}$, and then divide that by $\binom{N}{n}$ for drawing $n$
of the $N$ balls without regard to colour.  (All these choices are uniform.)

The expected value is just the same fraction of white balls in the sample as
in the urn:
\begin{equation}
   E[i] = n\frac{M}{N}
\end{equation}
and the variance is as follows:
\begin{equation}
   V[i] = n\frac{M(N-M)(N-n)}{N^2(N-1)}
\end{equation}
Proofs for mean and variance are in MathWorld~\cite{MathWorld}.

\section{Useful symmetries}

The Wikipedia article (as of this writing, of course; Wikipedia is a moving
target) gives some useful symmetries~\cite{Wikipedia}.  In our notation:
\begin{align}
   h(M,N,n,i) &= h(N-M,M,n,n-i) \label{first-symmetry} \\
   h(M,N,n,i) &= h(M,N,N-n,M-i) \label{second-symmetry} \\
   h(M,N,n,i) &= h(n,N,M,i) \label{third-symmetry}
\end{align}
If you have $M$ balls white, draw $n$, and hope for $i$ of them to be white,
you could instead flip all the colours, draw $n$, and hope for $n-i$ of them
to be white~\eqref{first-symmetry}.  Also, if you draw $n$ balls and find
$i$ to be white, that's the same as finding the $M-i$ remaining white balls
among the $N-n$ you did not draw; you can swap ``drawn'' and ``not drawn''
balls~\eqref{second-symmetry}.  Finally, you can swap the concepts of
``drawn'' and ``coloured white'' and imagine that the urn is choosing $M$
balls to possibly be drawn by you, instead of you choosing $n$ balls to
possibly be coloured white in the urn~\eqref{third-symmetry}.

\section{Tail inequalities}

We're interested in the chance that $i$ is at least $k$, for some $k$ that
will be a little bigger than the expected value $E[i]$.  We want to say that
when $k$ is just a tiny bit bigger than $E[i]$, then this chance is already
very small.  That will mean proving that this function is small:
\begin{equation}
   H(M,N,n,k) = \sum_{i=k}^{n} h(M,N,n,i) =
      \sum_{i=k}^{n} \binom{M}{i}\binom{N-M}{n-i}\left/\binom{N}{n}\right.
\end{equation}
That's the sum for all $i\ge k$ of the probability distribution
function~$h(M,N,n,i)$; we could equally correctly write the summation as
going to infinity, because $h(M,N,n,i)$ is zero for $i>n$; I wrote it up to
$n$ for consistency with Chv\'atal~\cite{Chvatal:Tail}.

Chv\'atal gives the following bound, which he credits to
Hoeffding~\cite{Chvatal:Tail,Hoeffding:Probability}.  I believe this is a
special case of the well-known result now known as Hoeffding's Inequality,
but that's a very powerful result and the steps required to apply it to the
hypergeometric distribution in particular are a little involved.  Where
$p=M/N$ and $k=(p+t)n$ with $t\ge 0$, we have this:
\begin{equation}
   H(M,N,n,k) \le \left( \left( \frac{p}{p+t} \right)^{p+t}
      \left( \frac{1-p}{1-p-t} \right)^{1-p-t} \right)^n
\end{equation}

That is a bit of a mess, but we can relax it a little further to get what
Chv\'atal describes as a ``more elegant but weaker'' bound which is more
likely what we'll want to use when applying this result:~\cite{Chvatal:Tail}
\begin{equation}
   H(M,N,n,k) \le e^{-2t^2n} \label{final-inequality}
\end{equation}

That's a nice one-sided tail inequality for hypergeometric distributions. 
Stating it in terms that sound like what we want for using it in proving
that a randomized algorithm works:  If $i$ is an hypergeometric random
variable with the parameters $N$, $M$, and $n$ as described above, then
\begin{equation}
   Pr[i\ge E[i]+tn] \le e^{-2t^2n}
\end{equation}

If we want an inequality for the other tail, then we can apply the
symmetry~\eqref{first-symmetry} as follows:
\begin{align}
   Pr[i\le k'] &= \sum_{i=0}^{k'} h(M,N,n,i) \\
      &= \sum_{i=0}^{k'} h(N-M,N,n,n-i)
\end{align}

Then we can change the index of summation to $j=n-i$ and get:
\begin{equation}
   Pr[i\le k'] = \sum_{j=n-k'}^{n} h(N-M,N,n,j)
\end{equation}

The other side's inequality~\eqref{final-inequality} can give us a nice
bound for that if we choose $k$, $t$, and $p$ properly.  We want
$k=n-k'=(p+t)n$ where $p=(N-M)/N=1-(M/N)$.  Then doing the algebra we get
$k'=E[i]-tn$, nicely equal and opposite to the other side's bound:
\begin{equation}
   Pr[i\le E[i]-tn] \le e^{-2t^2n}
\end{equation}

\bibliographystyle{plain}
\bibliography{hypergeometric}

\end{document}